# An asymptotic result for Brownian polymers

Thomas Mountford[a] and Pierre Tarrès[b,1]

[a]Ecole Polytechnique Fédérale de Lausanne, Département de mathématiques, 1015 Lausanne, Switzerland.
E-mail: thomas.mountford@epfl.ch
[b]University of Oxford, Mathematical Institute, 24-29 St Giles, Oxford OX1 3LB, United Kingdom.
E-mail: tarres@maths.ox.ac.uk

**Abstract.** We consider a model of the shape of a growing polymer introduced by Durrett and Rogers (*Probab. Theory Related Fields* **92** (1992) 337–349). We prove their conjecture about the asymptotic behavior of the underlying continuous process $X_t$ (corresponding to the location of the end of the polymer at time $t$) for a particular type of repelling interaction function without compact support.

**Résumé.** Nous considérons un modèle de formation de polymères introduit par Durrett et Rogers (*Probab. Theory Related Fields* **92** (1992) 337–349). Nous prouvons leur conjecture sur le comportement asymptotique du processus continu associé $X_t$ (correspondant à l'emplacement de l'extrémité du polymère au temps $t$) pour un type particulier de fonction d'interaction répulsive à support non compact.

## 1. Introduction

### 1.1. General setting

Let $(\Omega, \mathcal{F}, \mathbb{P})$ be a probability space, let $\{B_t, \mathcal{F}_t: t \geq 0\}$ be a Brownian motion on $\mathbb{R}^d$ (starting in 0 at time 0), and let $f: \mathbb{R}^d \to \mathbb{R}^d$ ($d \geq 1$) be a measurable function. In this paper we will consider processes $(X_t)_{t \geq 0}$ of the form

$$X_t = B_t + \int_0^t \mathrm{d}s \int_0^s f(X_s - X_u) \, \mathrm{d}u. \tag{1}$$

Equation (1) has a pathwise unique *strong* solution if $f$ is assumed to be Lipschitz, e.g. by Theorem 11.2 in [40]. Note that the assumption of $f$ being continuous is sufficient to conclude the uniqueness of a strong solution if it exists (see Theorem 5 and Corollary 1, p. 271 in [19]). The existence and uniqueness of a

---

[1]On leave from CNRS and Université Paul Sabatier, Toulouse, France

*weak* solution to (1) is ensured under the assumption that $f$ is bounded, using a generalization of Girsanov theorem (see Corollary 3.5.2 in [22]).

Observe that Eq. (1) is equivalent, in dimension one, to

$$X_t = B_t + \int_0^t \int_{-\infty}^{\infty} f(-z) L_s(X_s + z) \, dz \, ds \tag{2}$$

or

$$dX_t = dB_t + \left\{ \int_{-\infty}^{\infty} f(-z) L_t(X_t + z) \, dz \right\} dt,$$

where $L_t(y)$ is the *local time* (or *occupation time density*) of the process $X$. This formulation makes it clear how the process interacts with its own occupation density.

This setting has been introduced by Durrett and Rogers [18] in 1992 as a model for the shape of a growing polymer, $X_t$ corresponding to the location of the end of the polymer at time $t$. Without any assumption on the function $f$, the stochastic differential equation (1) defines a *self-interacting diffusion*, in the sense that the process $X$ evolves in an environment changing with its prior trajectory. We will call it *self-repelling* (resp. *self-attracting*) if, for all $x \in \mathbb{R}^d$, $x.f(x) \geq 0$ (resp. $\leq 0$), in other words if it is more likely to stay away from (resp. come back to) the places it has already visited before.

The model is a continuous analogue of the notion of edge (resp. vertex) self-interacting random walk (SIRW) on discrete graphs, defined as follows: at each step, the probability to move along an edge is proportional to a function – called the weight function – of the number of visits to this edge (resp. to the adjacent vertex). This notion was introduced in 1986 by Coppersmith and Diaconis [7] with the seminal definition of *edge-reinforced random walk* (ERRW) in the particular case of linear weight function.

Self-interacting random processes are useful in the understanding of self-organization and learning behaviour. Othmer and Stevens [33] suggest reinforced random walks as a model for the movement of myxobacteria (each bacterium moving to a site with a probability depending on the number of visits of all the bacteria to this site), and raise the question as to whether "aggregation is possible with such strictly local modification or whether some form of longer range communication is necessary". These reinforced walks can also describe spatial monopolistic competition in economics.

Let us mention some results on self-interacting diffusions (SIDs) and self-interacting random walks (SIRWs), and discuss the relationship between these two models.

*1.2. Previous results on self-interacting diffusions (SIDs)*

Durrett and Rogers [18] obtain upper bounds on the norm of the position of the particle for bounded functions $f$ of compact support and, in dimension one, lower bounds in the two following cases: $f$ nonnegative and $f$ bounded of the form $f(x) \sim lx^{-\beta}$ as $x \to \pm\infty$. Cranston and Mountford [9] describe precisely the asymptotic behaviour for nonnegative $f$, which solves a conjecture of Durrett and Rogers in [18] (the purpose of the present article is to prove their conjecture in the other case $f(x) \sim_{x \to \pm\infty} lx^{-\beta}$). The two articles [9, 18] are presented in more detail below.

The self-attracting case is studied in 1995 by Cranston and Le Jan [8]: the linear interaction (i.e. linear $f$) is considered, as well as the constant interaction (i.e. $f(x) = \sigma \text{sign}(x)$, $\sigma < 0$) in dimension one, both cases leading to an almost-sure convergence of the process, which is intuitive in the sense that self-attraction should lead to localization. The constant interaction result is generalized on $\mathbb{R}^d$ – $d \geq 2$ – (defined here by $f(x) = \sigma x/\|x\|$, $\sigma < 0$) by Raimond [38], again with a.s. convergence of the process. In dimension one, Herrmann and Roynette [20] generalize the a.s. convergence result to the case of odd, decreasing and bounded functions $f$ of lower bounded intensity in the neighbourhood of 0 (i.e. such that $|f(x)| \geq C\exp(-\rho/|x|^k)$ for small $|x|$, with constants $C, \rho > 0$ and $k \in \mathbb{N}$). The lower bounded intensity is indeed needed in order to always keep a strictly attractive force towards the already visited sites. In the case $f(x) = -\text{sign}(x)\mathbb{1}_{|x| \geq a}$, non-local in the sense that $f$ is zero in a neighbourhood of 0 so that the particle does not take into account

the visits to its close neighbourhood, Cranston and Le Jan [8] prove that the diffusion does not converge a.s., and Herrmann and Roynette [20] show that the paths are however a.s. bounded.

Let us also present alternative models of self-interacting diffusions. Norris, Rogers and Williams [32] define in 1987 a self-avoiding random walk as a Brownian motion model with local time drift. More recently, Benaïm, Ledoux and Raimond [2] introduce a model similar to Eq. (1), with this difference that the drift is given by an average of the past occupation (inserting a factor of $1/s$ in the first integral). Assuming that the particle lives in a compact connected smooth ($C^\infty$) Riemannian manifold (without boundary), and that $f(X_s - X_t)$ is replaced by the gradient of a potential $\nabla V_{X_s}(X_t)$ with sufficient differentiability, they prove that the normalized occupation measure $\mu_t = \frac{1}{t}\int_0^t \delta_{X_s}\,\mathrm{d}s$ asymptotically shadows the solutions of a deterministic differential equation, so that the possible limit sets of $\mu_t$ are "attractor free sets" for this equation. Depending on the structure of the interaction, various corresponding dynamics are possible; however, when the diffusion is self-repelling or weakly self-attracting, according to definitions introduced by the authors (taking into account that the particle lives in a compact set), $\mu_t$ a.s. converges toward the normalized Riemannian measure.

These self-interacting diffusions are further analysed by Benaïm and Raimond. In [3], convergence in law properties are discussed, whereas the symmetric interaction case is studied in [4]: $\mu_t$ converges almost surely to the critical set of a certain nonlinear free energy functional $J$. Generically, this critical set consists of finitely many points, so that $\mu_t$ converges a.s. toward a local minimum of $J$, each such minimum having a positive probability of being selected. A self-interacting model introduced by Del Moral and Miclo [13, 14] presents some similarity with the latter model: in a discrete time setting, the evolution depends on the present position and on the occupation measure created by the path up to this instant. The authors give sufficient conditions for a.s. convergence of the empirical measures, and provide upper bounds on the corresponding rate of convergence to the limiting measure.

### 1.3. Previous results on self-interacting random walks (SIRWs)

The results on discrete-time SIRWs are naturally ordered by the structure of the interaction: edge SIRW with respectively self-repelling, weakly reinforced, linear, once-reinforced or superlinear interaction, and vertex SIRW with linear and polynomial interaction.

The edge SIRW with self-repelling and weakly reinforced interactions, which correspond to decreasing and sublinear weight functions, have been studied by Tóth on the integer line $\mathbb{Z}$ in various cases, leading to results of convergence in law of the position of the random walk after renormalization: see [47] for the case of exponentially self-repelling random walks, which may be interpreted as the discrete-time counterpart of the self-repelling diffusion defined by (1), with $f$ odd of compact support; see also [48] for a survey of the different results obtained in this framework.

The linear edge SIRW (i.e. with linear weight function) corresponds to the critical case, and is generally called *edge-reinforced random walk* (ERRW). The ERRW on finite graphs is a mixture of reversible Markov chains, and the mixing measure can be determined explicitly ([15], see also [23, 41]), which has applications in Bayesian statistics [16]. On infinite graphs, the main question so far has been to give criteria for recurrence versus transience. Coppersmith and Diaconis [7] observe in 1986 that the walk is recurrent on $\mathbb{Z}$; more generally, on acyclic or directed graphs, the walk can be written as a random walk in an *independent* random environment, as was first observed by Pemantle in 1988, which enables to deduce a recurrence/transience phase transition on the binary tree [34] and recurrence/transience criteria or laws of large numbers in different instances [5, 24, 45]. In the case of infinite graphs with cycles, Merkl and Rolles have recently obtained recurrence criteria and asymptotic estimates on graphs of the form $\mathbb{Z} \times G$, $G$ finite graph, and a two-dimensional graph [27, 28, 30, 31, 42]. The fundamental question of recurrence or transience on $\mathbb{Z}^k$, $k \geq 2$, is still open.

The edge SIRW with one-time reinforcement, i.e. where the current weight of an edge is $1 + \delta$ if it has been crossed and 1 if it has never been crossed, is generally called *once-reinforced random walk* and has been introduced by Davis [10] in an attempt to provide a simplified version of the problem of establishing recurrence or transience for ERRW. The process is recurrent on ladders $\mathbb{Z} \times \{1,\ldots,d\}$ for $\delta \in (0, 1/(d-2))$ [43] and for large $\delta$ [50], and is transient on regular trees for all $\delta > 0$ or more generally on random trees generated by a supercritical branching process [6, 17], in contrast with the behavior of ERRW on these trees.

The superlinear edge SIRW has so far been studied under the condition of a reciprocally summable weight function. This condition is necessary and sufficient for visiting only a finite number of vertices in the case of nondecreasing weight functions. Davis [10] and Sellke [44] proved, respectively on the integer line and on $\mathbb{Z}^k$, $k \geq 2$, that this condition a.s. implies the existence of a random *attracting* edge. In the general case of a graph of bounded degree, the *attracting* edge property still holds for weight functions increasing like a power greater than 1, see Limic [25], and in fact for any reciprocally summable nondecreasing function, see Limic and Tarrès [26].

The vertex SIRW has so far mainly been considered in the case of a linear weight function, which is generally called *vertex-reinforced random walk* (VRRW) and was introduced by Pemantle in 1988 [35]. Vertex-reinforced random walks on finite complete graphs, with reinforcements weighted by factors associated to each edge of the graph, have been studied by Pemantle [36] and Benaïm [1]. On the integers $\mathbb{Z}$, Pemantle and Volkov showed that the VRRW a.s. visits only finitely many vertices and, with positive probability, eventually gets stuck on five vertices, and Tarrès [46] proved that this localization on five points is the almost sure behavior. On arbitrary graphs, Volkov [51] proved that VRRW localizes with positive probability on some specific finite subgraphs.

The vertex SIRW with weight function $W(n) = n^\rho$ has recently been studied by Volkov [52]. In the superlinear case $\rho > 1$, the walk a.s. visits two vertices infinitely often. In the sublinear case $\rho < 1$ the walk a.s. either visits infinitely many sites infinitely often or is transient; it is conjectured that the latter behaviour cannot occur, and that in fact all integers are infinitely often visited.

### 1.4. Link between the continuous and discrete cases

The relationship between the continuous and discrete cases has so far not been much investigated. Note that the interaction is local for SIRWs (in the sense that the particle only takes into account the visits to its neighbouring sites), whereas the SID generally evolves according to the past of the process on the whole space at any time, so that we a priori expect a more erratic or localized behaviour for SIRW, at least without renormalization in space.

Let us mention two continuous limits of SIRWs, for the exponentially self-repelling walk and the once-reinforced random walk. In the self-repelling case on $\mathbb{Z}$, Tóth and Werner [49] have constructed in 1998 a continuous process arising as limit of the renormalized position of the self-repelling walk with exponentially decaying weight function, called the *true self-repelling motion*. The self-repellance is "local" in the sense that it is only due to the occupation measure density at an "immediate neighbourhood", and the process is not solution of a stochastic differential equation, having finite variation of order $3/2$. Formally, this process is of the form (2) with $f := -\delta'$, the negative gradient of Dirac's delta, but without Brownian increments, the white noise disappearing in the scaling limit. This similarity with Eq. (2) led Tóth and Werner [49] to conjecture that self-interacting diffusions (1) with $f$ of compact support and exponentially self-repelling SIRWs display the same large scale asymptotic behaviour, the assumption on $f$ allowing the same local interaction mechanism.

As for once-reinforced random walks, Davis [11] establishes a connection between the weak limit of the walk and a diffusion which receives a push when at its maximum or minimum.

We are not aware of any study of self-interacting diffusions leading to results of the same nature.

Surveys on self-interacting random processes have been written by Davis [12], Merkl and Rolles [29], Pemantle [37] and Tóth [48], each viewing the subject from a different perspective.

### 1.5. Statement of the problem

From now on, we restrict ourselves to a process $(X_t)_{t \geq 0}$ taking values in $\mathbb{R}$, satisfying a stochastic differential equation (1), with $f : \mathbb{R} \to \mathbb{R}$ Lipschitz. Let us describe the results obtained by Durrett and Rogers [18] in 1992 in this one-dimensional setting. First, if $f$ is bounded and has a compact support, then there exists a constant $\Gamma < \infty$ so that

$$\limsup \frac{|X_t|}{t} \leq \Gamma \quad \text{a.s.}$$

One would like to say more about the existence of a limit for $X_t/t$. They prove that if $f$ is nonnegative and $f(0) > 0$, then there exists $\gamma > 0$ such that

$$\liminf \frac{X_t}{t} \geq \gamma.$$

Cranston and Mountford [9] have shown in 1996, under the weaker condition $f$ nonnegative and $f \not\equiv 0$, the strong law of large numbers for the polymer, i.e. that there exists a positive constant $c$ such that

$$\lim_{t \to \infty} \frac{X_t}{t} = c \quad \text{a.s.}$$

This result was partially conjectured in [18].

However, the assumption of $f$ being nonnegative is undesirable since it says that the process, neither repulsive nor attractive, keeps a steady drift towards the right and therefore always goes through new territory, so that its self-interaction will not modify its qualitative behaviour, but rather only its speed of convergence to infinity.

The situation where $f$ takes values of both signs and is "repulsive" ($\forall x \in \mathbb{R}$, $xf(x) \geq 0$) is more difficult to study since the particle, which avoids familiar territory, can in general receive contradictory signals from its left-hand and right-hand sides (respectively towards the right and the left), so that it does not necessarily decide on a direction in the long run. This repulsive case has so far led to a very small number of results. We first mention a conjecture of Durrett and Rogers in [18] in the case of an odd function $f$ of compact support.

**Conjecture (Durrett and Rogers [18]).** *Suppose $f$ has compact support, and $f(-x) = -f(x)$; then*

$$\frac{X_t}{t} \to 0 \quad a.s.$$

Tóth and Werner [49] also conjectured, by comparing this model with exponentially self-repelling random walks on $\mathbb{Z}$ [47], that under the same assumptions $X_t/t^{2/3}$ converges in law (see remarks in Section 1.4), which means that the particle has a super-diffusive behaviour despite the fact that it only looks at the time spent in its immediate neighbourhood.

When the function $f$ is not compactly supported, and in particular when $f$ is not integrable, the issue is different since the drift is expected to grow with time if the process keeps the same direction asymptotically. Let (A1), (A2) and (A3) be the following assumptions:

(A1) $|f(x)| \leq M$,
(A2) $f(x)$ is decreasing for $x \in [q, \infty)$,
(A3) $x^\beta f(x) \to l > 0$ as $x \to \infty$ with $0 < \beta < 1$.

Let us recall the heuristics described in the introduction of the article of Durrett and Rogers [18].
Letting $x_t = T^{-\alpha} X_{tT}$ and $W_t = T^{-1/2} B_{tT}$ we can rewrite (1) as

$$x_t = T^{1/2-\alpha} W_t + T^{2-\alpha} \int_0^t \mathrm{d}s \int_0^s f(T^\alpha (x_s - x_u)) \,\mathrm{d}u.$$

If we set

$$\alpha := \frac{2}{(1+\beta)}$$

so that $2 - \alpha = \alpha\beta$ and let $T \to \infty$ we expect that a possible limit (still called $x_t$ for simplicity) should satisfy

$$x_t = \int_0^t \mathrm{d}s \int_0^s \frac{l \, \mathrm{d}u}{(x_s - x_u)^\beta}.$$

One solution is $x_t = c_0 t^\alpha$ where $c_0$ satisfies

$$\alpha c_0^{\beta+1} = \int_0^1 \frac{l \, \mathrm{d}u}{(1 - u^\alpha)^\beta}. \tag{3}$$

Durrett and Rogers [18] obtained the two following Theorems A and B, and conjectured the following Theorem 1, which is the main result of this paper.

**Theorem A (Durrett and Rogers [18]).** *Suppose* (A1)–(A3) *hold and $\alpha$ and $c_0$ are as above. Then*

$$\limsup_{t \to \infty} \frac{X_t}{t^\alpha} \le c_0.$$

**Theorem B (Durrett and Rogers [18]).** *Suppose* (A1)–(A3) *hold, $f$ is nonnegative and $f(0) > 0$. Then*

$$\frac{X_t}{t^\alpha} \to c_0 \quad a.s.$$

**Theorem 1 (Conjecture 3 of Durrett and Rogers [18]).** *Suppose $f(x) = x/(1+|x|^{\beta+1})$ with $0 < \beta < 1$. Then with probability $1/2$,*

$$\frac{X_t}{t^\alpha} \to c_0.$$

The aim of this article is to prove Theorem 1. The conjecture is believable in that, given the long term nature of the function $f$, it is reasonable to think that once the motion has decided on a direction, the contribution to the drift in the opposing direction will become increasingly negligible.

We fix $0 < \beta < 1$ and assume that $f(x) = x/(1+|x|^{\beta+1})$ in the sequel.

## 2. Introduction to the ideas of the proof

### 2.1. Notation

Let $a \wedge b$ or $\min(a, b)$ (resp. $a \vee b$ or $\max(a, b)$) be the minimum (resp. the maximum) of $a$ and $b$. For all $x \in \mathbb{R}$, let $x^+ := \max(x, 0)$ and $x^- := \max(-x, 0)$.

Given $a$ and $b \in \mathbb{R}$ such that $b > a$, we use the convention that $[b, a]$ is the empty set.

Given a subset $A$ of $\mathbb{R}$, we let $A^* = A \setminus \{0\}$, $A_+ = \{x \in A \text{ and } x \ge 0\}$, $A_- = \{x \in A \text{ and } x \le 0\}$, $A_\pm^* = (A_\pm)^*$.

Let $\mathrm{Cst}(a_1, a_2, \ldots, a_p)$ denote a positive constant depending only on $a_1, a_2, \ldots a_p$, and let Cst denote a universal positive constant.

Given a real-valued function $g$ on $\mathbb{R}$, we let $\|g\|_\infty := \sup_{x \in \mathbb{R}} |g(x)|$.

We will make use of the following classical exponential inequality (cf. for instance [39], Proposition (1.8), p. 52):

$$\mathbb{P}\Big[\sup_{s \in [0,t]} B_s \ge a\Big] \le \exp\Big(-\frac{a^2}{2t}\Big). \tag{4}$$

*2.2. Sketch of the proof*

The first (and most important) step of the proof is the following proposition.

**Proposition 1.** $\mathbb{P}(\limsup_{t\to\infty} |X_t| = \infty) = 1$.

The proof of Proposition 1 relies on the particular shape of the drift when the process remains stuck in a bounded interval. More precisely, let us define, for any $u \in \mathbb{R}_+$ and any interval $I$,

$$h_u^I(x) = \int_0^u f(x - X_s) \mathbb{1}_{X_s \in I} \, ds, \qquad k_u^I(x) = \int_0^u f(x - X_s) \mathbb{1}_{X_s \notin I} \, ds,$$

$$g_u(x) = \int_0^u f(x - X_s) \, ds = h_u^I(x) + k_u^I(x).$$

Recall that

$$dX_u = dB_u + g_u(X_u) \, du.$$

Firstly, when the process remains in interval $I$ for a long time, $k_u^I$ remains constant, so that $h_u^I$ gives the main contribution to the drift. As long as $X_t$ does not leave $I$ its behavior is, on time intervals of fixed scale starting at a stopping time $S$, comparable to the behavior of a diffusion with drift $h_S^I$, as stated in the following simple Lemma 1.

Let us beforehand introduce the preliminary Definition 1.

**Definition 1.** *Let $I \subseteq \mathbb{R}$ be an interval, and let $S$ be an a.s. finite stopping time for filtration $(\mathcal{F}_t)_{t \geq 0}$.*
*We denote by $(B_t^S : t \geq 0)$ the Brownian motion $(B_{S+t} - B_S : t \geq 0)$ (which is independent of $\mathcal{F}_S$), by $(Z_t^S : t \geq 0)$ the time-shifted process $(X_{S+t} : t \geq 0)$, and by $(Y_t^S : t \geq 0)$ the diffusion*

$$Y_t^S = X_S + \int_0^t h_S^I(Y_u^S) \, du + B_t^S$$

*with a drift function $x \mapsto h_S^I(x)$ depending on the past, "frozen" at time $S$, of the process on interval $I$.*

**Lemma 1.** *Let $S$ be an a.s. finite stopping time for filtration $(\mathcal{F}_t)_{t \geq 0}$, let $I$ be an interval and let $v$ be a finite positive constant. Let $W_v(\mathbb{R})$ be the Wiener space of continuous paths $\omega : [0, v] \to \mathbb{R}$, equipped with the $\sigma$-algebra $\mathcal{G}$ generated by projection maps $\omega \mapsto \omega(t)$. Let $(Y_t^S)_{t \geq 0}$ and $(Z_t^S)_{t \geq 0}$ be the processes introduced in Definition 1.*
*Given $A \in \mathcal{G}$, assume that*

$$\mathbb{P}(Y_\cdot^S \in A | \mathcal{F}_S) \geq \varepsilon.$$

*Then*

$$\mathbb{P}(Z_\cdot^S \in A | \mathcal{F}_S) \geq \mathrm{Cst}(v, \varepsilon)$$

*a.s. on event $\{\|k_S^I\|_\infty \leq v\}$.*

Lemma 1 is an easy consequence of Girsanov theorem, and is proved in Section 3.2.

Secondly, given a fixed time $u$ and a space interval $I = [a, b]$, the drift function $x \mapsto h_u^I(x)$ satisfies the following property: when $x \in \mathbb{R}$ is close to the right boundary of $I$, according to a definition involving only $a$ and $b$, then either $h_u^I(x)$ is positive or $h_u^I(y)$ is nonpositive for all $y \in [a, x]$, as implied by the following Lemma 2.

Let $x_{\max} := (1/\beta)^{1/(1+\beta)}$ be the point of change of monotonicity of $f$.

**Lemma 2.** *Let $u \in \mathbb{R}_+^*$, $a$, $b \in \mathbb{R}$. Suppose there exists $x_0 \in [a, b]$ such that $h_u^{[a,b]}(x_0) \leq 0$, and either $f(b - x_0) \leq f(b - a)^2$ and $b - x_0 \leq 1/16$, or $b - a \leq x_{\max}$. Then, for all $x \in [a, x_0]$, $h_u^{[a,b]}(x) \leq 0$.*

Lemma 2 is proved in Section 3.3. The symmetrical statement at the left boundary holds similarly, replacing $\leq$ by $\geq$ in the inequalities involving $h_u^{[a,b]}(x_0)$ and $h_u^{[a,b]}(x)$, and $b - x_0$ by $x_0 - a$. An equivalent statement to Lemma 2 is that there exists a constant $c \in (a, b)$ only depending on $a$, $b$ and $\beta$ such that, either $h_u(x) \geq 0$ for all $x \in [c, b]$, or there exists $x_0 \in [c, b]$ such that $h_u(x)(x - x_0) \geq 0$ for all $x \in [a, b]$. A consequence of Lemmas 1 and 2 is that, as long as $X_t$ remains in $I$, each time it approaches the border of $I$ the probability to leave it within a time limit depending only on the size of the interval is lower bounded. Therefore, the range of the process $X_t$ regularly widens, which explains Proposition 1, proved in Section 3.4.

In the remainder of this introduction, let us assume that $\{\limsup X_t = \infty\}$ holds. The treatment of the event $\{\liminf X_t = -\infty\}$ is similar, by symmetry.

Next, we prove that each time $X_t$ reaches its maximum, the probability that it surpasses its supremum by one within one unit of time is lower bounded, independently of the prior occupation measure of the process. This observation leads to the following lemma. For all $x \in \mathbb{R}$, let the stopping time $T_x$ be the first hitting time of $x$ by the process $(X_t : t \geq 0)$. For all $x \in \mathbb{R}_+$, let us define the following event

$$A_x := \{T_x \leq T_{x-1} + 1 < \infty\}.$$

Then $A_x$ holds for infinitely many $x \in \mathbb{Z}_+$, as stated in Lemma 3.

**Lemma 3.** *One has*

$$\{\limsup X_t = \infty\} \subseteq \left\{\sum_{x \in \mathbb{Z}_+} \mathbb{1}_{A_x} = \infty\right\}.$$

Lemma 3 is proved in Section 4.1.
Let us now use notation

$$G(t) := g_t(X_t)$$

for the drift at time $t$.

The introduction of events $A_x$ in the study of $X_t$ is justified by the following observation: For sufficiently large $x \in \mathbb{Z}_+$, if $A_x$ holds then the drift at time $T_x$ can be lower bounded by a term of the order of $T_x^{\alpha-1}$ using that, by Theorem A, $|X_t|/t^\alpha$ has remained upper bounded by $2c_0$ for $t \in [t_0, T_x]$, $t_0$ sufficiently large. This is the purpose of the following Lemma 4.

**Lemma 4.** *One has*

$$\liminf_{x \to \infty, T_x \leq T_{x-1}+1} \frac{G(T_x)}{T_x^{\alpha-1}} \geq \liminf_{x \to \infty, T_x \leq T_{x-1}+1} \frac{\inf_{y \in [x-1/2, x+1]} g_{T_x}(y)}{T_x^{\alpha-1}} \geq (4c_0)^{-\beta} > 0 \quad a.s.$$

**Proof.** The first inequality is straightforward.

Let us prove the second inequality. Let $x \geq 1$ be such that $T_x \leq T_{x-1} + 1$. Let us first make use of Theorem A: Fix $t_0 \in \mathbb{Z}_+$ such that, for all $t \in \mathbb{R}_+$, $|X_t| \leq 2c_0(t \vee t_0)^\alpha$. Note that this implies

$$x - \inf_{0 \leq s \leq T_x} X_s = X_{T_x} - \inf_{0 \leq s \leq T_x} X_s \leq 4c_0(T_x \vee t_0)^\alpha. \tag{5}$$

Suppose $A_x$ holds, and let $y \in [x - 1/2, x + 1]$. Then, using Eq. (5), $\|f\|_\infty \leq 1$ and $f(\gamma) \geq f(\alpha) \wedge f(\theta)$ for all $0 \leq \alpha \leq \gamma \leq \theta$ (see Section 3.1), we deduce that

$$g_{T_x}(y) = \int_0^{T_x} f(y - X_s) \, ds \geq \int_0^{T_{x-1}} f(y - X_s) \, ds - 1$$

$$\geq \Big(\inf_{z \in [\inf_{s \leq t} X_s, x-1]} f(y-z)\Big) T_{x-1} - 1$$

$$\geq \min\Big(f\Big(\frac{1}{2}\Big), f(4c_0(T_x \wedge t_0)^\alpha + 1)\Big)(T_x - 1) - 1.$$

This enables us to conclude, using that $f(u)/u^{-\beta} \to_{u \to \infty} 1$ and $2 - \alpha = \alpha\beta$. □

The last step makes use of the proof of Theorem 4 in [18] (stated here as Theorem B). This theorem cannot be applied directly, since the assumption $f$ nonnegative and $f(0) > 0$ is not satisfied. However, Lemma 4 will imply that the drift increases fast asymptotically and that event $A_x$ occurs for any large $x$, so that the negative contributions $\int_0^t (f(X_t - X_s))^- \, ds$ in $G(t)$ become increasingly negligible. More precisely, let

$$\Gamma_+ := \Big\{\liminf_{t \to \infty} \Big(\inf_{0 \leq s \leq t} (g_t(X_t) - g_s(X_t))\Big) > -\infty\Big\} \cap \{\liminf G(t) > 0\}.$$

Then $\Gamma_+$ a.s. holds on $\{\limsup X_t = \infty\}$, as implied by the following Lemma 5, and the proof of Theorem 4 in [18] can be adapted to show that $X_t/t^\alpha$ a.s. converges to $c_0$ on $\Gamma_+$, as stated hereafter in Lemma 6. The two Lemmas 5 and 6 complete the proof of Theorem 1.

**Lemma 5.** *One has*

$$\{\limsup X_t = \infty\} \subseteq \Gamma_+ \cap \Big\{\liminf \frac{G(t)}{t^{\alpha-1}} > 0\Big\} \quad a.s.$$

**Lemma 6.** *One has*

$$\Gamma_+ \subseteq \Big\{\lim \frac{X_t}{t^\alpha} = c_0\Big\} \quad a.s.$$

Lemma 5 is proved in Section 4.2; it provides a tighter asymptotic estimate of $G(t)$ than the one occurring on $\Gamma_+$ but this improvement, which is a consequence of the proof, is not required in Lemma 6. Lemma 6 is proved in Section 4.3.

*2.3. Outline of contents*

Section 3 begins with some remarks on function $f$ in Section 3.1, and then provides the proofs of Lemmas 1, 2 and Proposition 1, in Sections 3.2–3.4. Section 4 is devoted to the proofs of Lemmas 3, 5 and 6, in Sections 4.1–4.3.

## 3. Proof of Lemmas 1, 2 and Proposition 1

*3.1. Some remarks on function f*

We will need in the proof some observations about the function

$$f(x) = \frac{x}{1 + |x|^{1+\beta}}.$$

First, for $x \in \mathbb{R}_+$, the derivative of $f$ is

$$f'(x) = \frac{1 - \beta x^{1+\beta}}{(1 + x^{1+\beta})^2}.$$

- Hence, on $\mathbb{R}_+$, $f$ increases until $x_{\max} := (1/\beta)^{1/(1+\beta)} \geq 1$, and decreases after $x_{\max}$; remark that $f(x_{\max}) \leq x_{\max}^{-\beta} \leq 1$. Therefore $\|f\|_\infty \leq 1$.
- For all $\alpha, \gamma, \theta \in [-x_{\max}, \infty)$ such that $\alpha \leq \gamma \leq \theta$,

$$f(\gamma) \geq \min(f(\alpha), f(\theta)). \tag{6}$$

- Let us prove that $f'(x) \geq -f(x)$ for all $x \in [-1/2, \infty)$.

The inequality is straightforward on $[0, 1]$, since $f$ and $f'$ are both nonnegative on this interval. Observe that, for all $x \geq 1$,

$$\left|\frac{f'(x)}{f(x)}\right| = \left|\frac{1 - \beta x^{1+\beta}}{x(1 + x^{1+\beta})}\right| \leq \frac{1 + \beta x^{1+\beta}}{1 + x^{1+\beta}} \leq 1.$$

It remains to study the case $x \in [-1/2, 0]$. Since $f$ is odd, it suffices to prove that $f'(x) \geq f(x)$ for all $x \in [0, 1/2]$. On $[0, 1/2]$, $f'$ is decreasing and $f$ is increasing, and therefore $f'/f$ is decreasing. Now

$$\frac{f'(1/2)}{f(1/2)} = \frac{2(1 - (\beta/2^\beta)(1/2))}{1 + 1/2^{(1+\beta)}} \geq \frac{2(1 - 1/4)}{1 + 1/2} = 1,$$

using $\beta/2^\beta \leq 1/2$ (since $f(x) = x/2^x$ increases on $[0, 1]$). This yields the result.

### 3.2. Proof of Lemma 1

Assume that $\mathbb{P}(Y_\cdot^S \in A | \mathcal{F}_S) \geq \varepsilon$. One can write the process $Z_\cdot^S$ as

$$Z_t^S = X_S + \int_0^t g_{S+u}(Z_u^S)\,\mathrm{d}u + B_t^S$$

$$= X_S + \int_0^t h_S^I(Z_u^S)\,\mathrm{d}u + U_t^S,$$

where $U_t^S := B_t^S + \int_0^t C_t^S\,\mathrm{d}t$, with $C_u^S := g_{S+u}(Z_u^S) - h_S^I(Z_u^S)$. Note that $C_u^S$ is $\mathcal{F}_{S+u}$-measurable.

Conditioned on $\mathcal{F}_S$, the solutions $Y_\cdot^S$, $Z_\cdot^S$ are strong, adapted to the filtration $(\mathcal{F}_{S+t})_{t \geq 0}$ and are a.s. equal to a Borel function of respectively the pair $(X_S, B_\cdot^S)$ and $(X_S, U_\cdot^S)$ arising out of the limit of Picard iteration. See e.g. [39], Chapter IX. Thus the events $\{Y_\cdot^S \in A\}$ and $\{Z_\cdot^S \in A\}$ can be rewritten as respectively $\{(X_S, B_\cdot^S) \in A'\}$ and $\{(X_S, U_\cdot^S) \in A'\}$.

By Girsanov's theorem (see e.g. [39], Chapter VIII), there exists a probability $\widetilde{\mathbb{P}}$ under which the process $(U_t^S : 0 \leq t \leq v)$ has Wiener measure as its distribution (so that

$$\widetilde{\mathbb{P}}((X_S, U_\cdot^S) \in A' | \mathcal{F}_S) = \mathbb{P}((X_S, B_\cdot^S) \in A' | \mathcal{F}_S) \tag{7}$$

in particular) and

$$\mathbb{P}(Z_\cdot^S \in A | \mathcal{F}_S) = \mathbb{P}((X_S, U_\cdot^S) \in A' | \mathcal{F}_S) = \widetilde{\mathbb{E}}(\mathbb{1}_{\{(X_S, U_\cdot^S) \in A'\}} \exp(\Delta) | \mathcal{F}_S), \tag{8}$$

where

$$\Delta := \int_0^v C_u^S\,\mathrm{d}U_u^S - \frac{1}{2}\int_0^v (C_u^S)^2\,\mathrm{d}u.$$

Let us assume that $\{\|k_S^I\|_\infty \leq v\}$ holds. Our goal is to estimate from below $\Delta$ on an event of large probability. To this end, observe that, for $u \in [0, v]$,

$$|C_u^S| \leq |k_S^I(Z_u^S)| + |g_{S+u}(Z_u^S) - g_S(Z_u^S)| \leq \|k_S^I\|_\infty + v\|f\|_\infty \leq 2v,$$

using that $\|f\|_\infty \leq 1$.

Now, the process $(M_t^S)_{0 \leq t \leq v}$ defined by

$$M_t^S := \exp\left(-\int_0^t C_u^S \, \mathrm{d}U_u^S - \frac{1}{2}\int_0^t (C_u^S)^2 \, \mathrm{d}u\right)$$

is a martingale on $(\Omega, (\mathcal{F}_{S+t})_{0 \leq t \leq v}, \widetilde{\mathbb{P}})$, and therefore

$$\widetilde{\mathbb{P}}(\Delta \leq -x|\mathcal{F}_S) \leq \widetilde{\mathbb{P}}(M_v^S \geq \mathrm{e}^{x-4v^3}|\mathcal{F}_S) \leq \widetilde{\mathbb{E}}(M_v^S|\mathcal{F}_S)\mathrm{e}^{4v^3-x} = \mathrm{e}^{4v^3-x} = \frac{\varepsilon}{2}, \tag{9}$$

choosing $x := 4v^3 - \ln(\varepsilon/2)$. Consequently, Eqs (7)–(9) imply (recall that $\mathbb{P}(Y_\cdot^S \in A|\mathcal{F}_S) \geq \varepsilon$ by assumption)

$$\begin{aligned}
\mathbb{P}(Z_\cdot^S \in A|\mathcal{F}_S) &\geq \mathrm{e}^{-x}\widetilde{\mathbb{E}}[\mathbb{1}_{\{(X_S, U_\cdot^S) \in A'\}}\mathbb{1}_{\{\Delta > -x\}}|\mathcal{F}_S] \\
&\geq \mathrm{e}^{-x}[\widetilde{\mathbb{P}}((X_S, U_\cdot^S) \in A'|\mathcal{F}_S) - \widetilde{\mathbb{P}}(\Delta \leq -x|\mathcal{F}_S)] \\
&\geq \mathrm{e}^{-x}\left[\mathbb{P}((X_S, B_\cdot^S) \in A'|\mathcal{F}_S) - \frac{\varepsilon}{2}\right] = \mathrm{e}^{-x}\left[\mathbb{P}(Y_\cdot^S \in A|\mathcal{F}_S) - \frac{\varepsilon}{2}\right] \\
&\geq \mathrm{e}^{-x}\frac{\varepsilon}{2} = \mathrm{e}^{-4v^3}\frac{\varepsilon^2}{4}.
\end{aligned}$$

### 3.3. Proof of Lemma 2

Let $h_u := h_u^{[a,b]}$ for simplicity.

Let us assume that $b - a > x_{\max}$: indeed, when $b - a \leq x_{\max}$, $h_u$ is nondecreasing on $[a, b]$, since $f$ is increasing on $[-x_{\max}, x_{\max}]$, and hence $h_u(x) \leq h_u(x_0) \leq 0$ for all $x \in [a, x_0]$.

Let $x_0 \in [b - 1/16, b]$ be such that $h_u(x_0) \leq 0$ and $f(b - x_0) \leq f(b - a)^2$.

Firstly, let us prove that, for all $x \in [b - 1/2, x_0]$, $h_u'(x) \geq -h_u(x)$. This will imply that $x \mapsto h_u(x)\mathrm{e}^x$ increases on $[b - 1/2, x_0]$, and therefore that $h_u(x) \leq 0$ for all $x \in [b - 1/2, x_0]$.

Indeed, let $x \in [b - 1/2, x_0]$. Then

$$h_u'(x) = \int_0^u f'(x - X_s)\mathbb{1}_{X_s \in [a,b]} \, \mathrm{d}s.$$

Now, if $X_s \in [a, b]$ then $x - X_s \in [-1/2, \infty)$ and therefore $f'(x - X_s) \geq -f(x - X_s)$ (see the remarks on function $f$, Section 3.1). Hence

$$h_u'(x) \geq -\int_0^u f(x - X_s)\mathbb{1}_{X_s \in [a,b]} \, \mathrm{d}s = -h_u(x).$$

Let us now consider the case $x \in [a, b - 1/2]$. Let us compare $h_u(x)$ to $h_u(x_0)$. Observe that

$$\begin{aligned}
h_u(x) &= \int_0^u f(x - X_s)\mathbb{1}_{X_s \in [a,x] \cup [x_0,b]} \, \mathrm{d}s + \int_0^u f(x - X_s)\mathbb{1}_{X_s \in (x,x_0)} \, \mathrm{d}s \\
&\leq \int_0^u f(x - X_s)\mathbb{1}_{X_s \in [a,x]} \, \mathrm{d}s - \int_0^u f(X_s - x)\mathbb{1}_{X_s \in [x_0,b]} \, \mathrm{d}s \tag{10}
\end{aligned}$$

since $f(x - X_s) \leq 0$ when $X_s \in [x, x_0]$.

On the other hand, using a similar argument,

$$0 \geq h_u(x_0) \geq \int_0^u f(x_0 - X_s)\mathbb{1}_{X_s \in [a,x]} \, \mathrm{d}s - \int_0^u f(X_s - x_0)\mathbb{1}_{X_s \in [x_0,b]} \, \mathrm{d}s. \tag{11}$$

Now, $x_0 - x \geq 7/16$ (since $x \leq b - 1/2$ and $x_0 \geq b - 1/16$) and $x_0$, $x \in [a,b]$ imply that, when $X_s \in [a,x]$, then $x_0 - X_s \in [7/16, b-a]$ and that, when $X_s \in [x_0, b]$, $X_s - x_0 \in [0, b-x_0] \subseteq [0, x_{\max}]$ and $X_s - x \in [7/16, b-a]$.

Therefore, Eqs (10) and (11) yield respectively, together with $\|f\|_\infty \leq 1$ and inequality (6),

$$h_u(x) \leq \int_0^u \mathbb{1}_{X_s \in [a,x]} \, ds - \min\left(f\left(\frac{7}{16}\right), f(b-a)\right) \int_0^u \mathbb{1}_{X_s \in [x_0,b]} \, ds$$

and

$$0 \geq h_u(x_0) \geq \min\left(f\left(\frac{7}{16}\right), f(b-a)\right) \int_0^u \mathbb{1}_{X_s \in [a,x]} \, ds - f(b-x_0) \int_0^u \mathbb{1}_{X_s \in [x_0,b]} \, ds.$$

Thus, it suffices to show that $(\min(f(7/16), f(b-a)))^2 \geq f(b-x_0)$ to conclude that $h_u(x) \leq 0$. We already know that $f(b-a)^2 \geq f(b-x_0)$ by assumption, and thus it remains to prove that $f(7/16)^2 \geq f(1/16)$, since $f(1/16) \geq f(b-x_0)$ (recall that $b - x_0 \in [0, 1/16]$ by assumption). This inequality follows from

$$f\left(\frac{7}{16}\right) \geq f\left(\frac{3}{8}\right) \geq \frac{3}{8} \frac{1}{(1+1/2^{1+\beta})} \geq \frac{3}{8} \frac{1}{(3/2)} = \frac{1}{4},$$

and $f(1/16) \leq 1/16$.

### 3.4. Proof of Proposition 1

Let us define the sequence $(a_n)_{n \in \mathbb{Z}_+}$ by $a_0 := 0$, $a_1 := x_{\max}/2$ and, recursively, for all $n \geq 1$,

$$a_{n+1} := a_n + \frac{1}{2}\left(f(4a_{n-1})^2 \wedge \frac{1}{16}\right).$$

It is immediate from the definition that $(a_n)_{n \in \mathbb{Z}}$ is an increasing sequence, converging to infinity as $n \to \infty$, and such that for all $n \geq 1$, $a_{n+1} \leq 2a_n$ (using that $a_n \geq x_{\max}/2 \geq 1/2$).

For all $n \in \mathbb{Z}_+$ and $t \in \mathbb{R}_+$, let us define the stopping time

$$S_{n,t} := \inf\{u > t \text{ s.t. } X_u \notin [-a_n, a_n]\};$$

remark that $S_{0,t} = t$.

Let us prove by induction on $n \in \mathbb{Z}_+$ that, for all $n \in \mathbb{Z}_+$ and $t \in \mathbb{R}_+$, $S_{n,t} < \infty$ a.s. This will obviously suffice to prove the proposition. The case $n = 0$ is trivial.

Let $n \geq 1$ and $t \in \mathbb{R}_+$. We will prove that there exists a positive constant $\zeta_{n,t}$ (depending only on $n$ and $t$) such that, for all $s \geq t$,

$$\mathbb{P}(S_{n,t} < \infty | \mathcal{F}_{S_{n-1,s}}) \geq \zeta_{n,t} > 0 \quad \text{a.s. on } \{S_{n-1,s} < \infty\}. \tag{12}$$

This will enable us to complete the induction step: indeed, for all $t \in \mathbb{R}_+$ and $s \geq t$, the induction assumption implies $S_{n-1,s} < \infty$ a.s., and we deduce from (12) that $\mathbb{E}(\mathbb{1}_{S_{n,t}<\infty}|\mathcal{F}_s) \geq \zeta_{n,t}$ a.s. By a standard martingale convergence theorem, $\mathbb{E}(\mathbb{1}_{S_{n,t}<\infty}|\mathcal{F}_s) \to_{s \to \infty} \mathbb{1}_{S_{n,t}<\infty} \geq \zeta_{n,t} > 0$, and therefore $S_{n,t} < \infty$ a.s.

Given $t \in \mathbb{R}_+$ and $s \geq t$, let us now prove (12). Let $S := S_{n-1,s}$, $I := [-a_n, a_n]$ and $h_S := h_S^I$. If $S \geq S_{n,t}$, there is nothing to prove. Hence we assume $S < S_{n,t}$, which implies $X_S \in [-a_n, -a_{n-1}] \cup [a_{n-1}, a_n]$. We do the proof under the assumption $X_S \in [a_{n-1}, a_n]$, the treatment of the case $X_S \in [-a_n, -a_{n-1}]$ being similar.

Let $(B_t^S)_{t \geq 0}$ and $(Y_t^S)_{t \geq 0}$ be the processes introduced in Definition 1 (with $S$ and $I$ fixed above).

Consider the stopping time

$$\tau := \inf\{u > 0 \text{ s.t. } Y_u^S \notin [-a_n, a_n]\}.$$

We now prove that

$$\mathbb{P}(\tau \leq 8a_n^2 \vee 1 |\ \mathcal{F}_S) \geq \frac{1}{2}, \tag{13}$$

which implies by Lemma 1 that $\mathbb{P}(S_{n,t} \leq S + 8a_n^2 \vee 1|\mathcal{F}_S) \geq \mathrm{Cst}(a_n, t) > 0$ a.s., using that $\|k_S^{[-a_n,a_n]}\|_\infty \leq t$ (as a consequence of $S < S_{n,t}$ and $\|f\|_\infty \leq 1$), and therefore yields (12).

We distinguish between two cases (note that $h_S(a_n) > 0$, and that $2a_{n-1} - a_n = a_{n-1} - (a_n - a_{n-1})$):

(1) $\forall x \in [2a_{n-1} - a_n, a_n],\ h_S(x) \geq 0$.

Then the process $(Y_u^S)_{u \geq 0}$ is, before leaving $[2a_{n-1} - a_n, a_n]$, bigger than the Brownian motion $(B_u^S)_{u \geq 0}$ plus $a_{n-1}$, and therefore the probability for $Y_u^S$ to leave the interval $[2a_{n-1} - a_n, a_n]$ at $a_n$ before time 1 is greater than or equal to the probability of $B_u^S$ leaving the interval $[-(a_n - a_{n-1}), a_n - a_{n-1}]$ at $a_n - a_{n-1}$ before time 1 (using a comparison result, cf. [21] for instance). Hence, using $a_n - a_{n-1} \leq 1$,

$$\mathbb{P}(\tau \leq 1|\mathcal{F}_S) \geq \frac{1}{2}\mathbb{P}(|B_1^S| \geq a_n - a_{n-1}|\mathcal{F}_S) \geq \frac{1}{2}\mathbb{P}(|B_1^S| \geq 1|\mathcal{F}_S) = \mathrm{Cst}. \tag{14}$$

(2) $\exists x \in [2a_{n-1} - a_n, a_n]$ such that $h_S(x) = 0$.

Let $x_0$ be the greatest $x$ in $[2a_{n-1} - a_n, a_n]$ such that $h_S(x) = 0$. Next apply Lemma 2, with $a := -a_n$ and $b := a_n$. Let us check that the assumptions are fulfilled. If $n = 1$, then $b - a = 2a_n \leq x_{\max}$. If $n \geq 2$, then it follows from the recursive definition of $(a_n)_{n \in \mathbb{Z}_+}$ that $b - x_0 \leq 2(a_n - a_{n-1}) \leq 1/16$, and

$$f(b - x_0) \leq b - x_0 \leq 2(a_n - a_{n-1}) \leq f(4a_{n-1})^2 \leq f(2a_n)^2,$$

using in the last inequality that $x_{\max} \leq 2a_n \leq 4a_{n-1}$ and that $f$ is decreasing on $[x_{\max}, \infty)$.

Hence, Lemma 2 implies $h_S(x) \leq 0$ for all $x \in [-a_n, x_0]$. Consequently, using that $h_S(x) \geq 0$ for all $x \in [x_0, a_n]$ by definition of $x_0$, the process $(R_u^S)_{u \geq 0}$, defined by

$$R_u^S = (Y_{u \wedge \tau}^S - x_0)^2 - u \wedge \tau$$

is a submartingale. Therefore

$$0 \leq \mathbb{E}(R_{8a_n^2}^S - R_0^S|\mathcal{F}_S) \leq (2a_n)^2 - 8a_n^2 \mathbb{P}(\tau \geq 8a_n^2),$$

thus

$$\mathbb{P}(\tau \geq 8a_n^2|\mathcal{F}_S) \leq \frac{1}{2}. \tag{15}$$

## 4. Proof of Lemmas 3, 5 and 6

### 4.1. Proof of Lemma 3

Let us first introduce Definition 2 and Lemma 7.

**Definition 2.** *For all $a, t \in \mathbb{R}_+$ and any continuous process $(R_t)_{t \geq 0}$, let $U(a, t, R)$ be the stopping time*

$$U(a, t, R) := \inf\{u > t\ s.t.\ |R_u - R_t| = a\} \wedge \left(t + \frac{a^2}{2}\right).$$

**Lemma 7.** *Let $h \in \mathbb{R}_+$, and let $M_t := B_t + ht$ be a Brownian motion with drift. Then, for every a.s. finite stopping time $S$ and all $a \in \mathbb{R}_+^*$ such that $ah \geq 6$,*

$$\mathbb{P}[M_{U(a,S,M)} = M_S + a|\mathcal{F}_t] \geq 1 - \exp(2 - ah).$$

**Proof.** Assume that $ah \geq 6$, and that $S = 0$ for simplicity. Let $U := U(a, 0, M)$.

For all $\lambda \in \mathbb{R}$, $N_t^\lambda := \exp(\lambda B_t - \lambda^2 t/2) = \exp(\lambda M_t - \lambda(\lambda/2 + h)t)$ is a martingale. We choose $\lambda := 2/a - h$, which satisfies $\lambda < 0$ and $\lambda/2 + h > 0$. Then

$$1 = \mathbb{E}[N_U^\lambda | \mathcal{F}_0] \geq \mathbb{P}\left[M_U = -a \text{ or } U = \frac{a^2}{2} \Big| \mathcal{F}_S\right] \left[\exp(-\lambda a) \wedge \exp\left(\lambda a - \lambda\left(\frac{\lambda}{2} + h\right)\frac{a^2}{2}\right)\right]$$

$$= \mathbb{P}\left[M_U = -a \text{ or } U = \frac{a^2}{2} \Big| \mathcal{F}_S\right] \left[\exp(ah - 2) \wedge \exp\left(\frac{(ah-2)^2}{4}\right)\right]$$

$$= \mathbb{P}\left[M_U = -a \text{ or } U = \frac{a^2}{2} \Big| \mathcal{F}_S\right] \exp(ah - 2).$$

□

Let

$$\mathcal{E} := \{\limsup X_t < \infty\} \cup \left\{\sum_{x \in \mathbb{Z}} \mathbb{1}_{A_x} = \infty\right\}$$

and, for all $y \in \mathbb{Z}_+$,

$$\mathcal{E}_y := \{\limsup X_t < \infty\} \cup \left\{\sum_{x \geq y} \mathbb{1}_{A_x} \geq 1\right\}.$$

Then $\mathcal{E} = \bigcap_{y \in \mathbb{Z}_+} \mathcal{E}_y$ and, for all $z \geq y$, $\mathcal{E}_z \subseteq \mathcal{E}_y$.

Our goal is to prove that, for all $y \in \mathbb{Z}_+$,

$$\mathbb{P}(\mathcal{E}_y | \mathcal{F}_{T_y}) \geq \text{Cst} > 0 \quad \text{a.s. on } \{T_y < \infty\}. \tag{16}$$

This enables us to conclude that $\mathbb{P}(\mathcal{E}) = 1$. Indeed, this implies that for all $y \in \mathbb{Z}_+$ and $s \in \mathbb{R}_+$, $\mathbb{P}(\mathcal{E}_y | \mathcal{F}_s) \geq$ Cst $> 0$ since, if $z := \inf\{n \in \mathbb{Z}_+ \text{ s.t. } n \geq \sup_{0 \leq t \leq s} X_t \vee y\}$, then

$$\mathbb{P}(\mathcal{E}_y | \mathcal{F}_s) \geq \mathbb{P}(\mathcal{E}_z | \mathcal{F}_s) \geq \mathbb{E}(\mathbb{E}(\mathbb{1}_{\mathcal{E}_z} | \mathcal{F}_{T_z}) \mathbb{1}_{T_z < \infty} | \mathcal{F}_s) + \mathbb{P}(T_z = \infty | \mathcal{F}_s) \geq \text{Cst}.$$

We deduce subsequently that, for all $y \in \mathbb{Z}_+$, $\mathcal{E}_y$ holds a.s. (using that $\mathbb{P}(\mathcal{E}_y | \mathcal{F}_s)$ tends to $\mathbb{1}_{\mathcal{E}_y}$ a.s. as $s$ tends to infinity), which completes the proof of $\mathbb{P}(\mathcal{E}) = 1$.

Given $y \in \mathbb{Z}_+$, let us prove Eq. (16). Almost surely on $\{T_y < \infty\}$, the distribution of $X_t$, $0 \leq t \leq T_y$, is absolutely continuous with respect to that of Brownian motion run until hitting $y$, and therefore has a.s. a bounded local time $l_y$.

Assume that $T_y < \infty$. Let us define

$$\delta := \inf\left\{u > 0 \text{ s.t. } \int_{y-u}^{y} l_y(v)\, dv = \frac{1}{u^2}\right\} \wedge 1.$$

Note that $\delta$ exists, since $\lim_{u \to \infty} 1/u^2 = 0$ and $\int_{-\infty}^{y} l_y(v)\, dv = T_y$.

Let $U(a, t, X)$ be the stopping time of Definition 2. Observe that

$$\bigcap_{n \in \mathbb{Z}_+, n\delta < 1} \{T_{y+(n+1)\delta} = U(\delta, T_{y+n\delta}, X)\}$$

$$\subseteq \left\{T_{y+1} \leq T_y + \left(\sup_{k \in \mathbb{Z}_+, k\delta < 1} k + 1\right)\frac{\delta^2}{2} \leq T_y + 1\right\} \subseteq \mathcal{E}_y. \tag{17}$$

Therefore, it suffices to estimate from below, for all $n \in \mathbb{Z}_+$ such that $n\delta < 1$,

$$\mathbb{P}(T_{y+(n+1)\delta} = U(\delta, T_{y+n\delta}, X)|\mathcal{F}_{T_{y+n\delta}})$$

provided that $\Delta_n := \bigcap_{0 \leq k \leq n-1}\{T_{y+(k+1)\delta} = U(\delta, T_{y+k\delta}, X)\}$ holds (which implies $T_{y+n\delta} < \infty$). To achieve this goal, we give a lower bound of the drift $g_t(X_t)$ on the time interval $[T_{y+n\delta}, U(\delta, T_{y+n\delta}, X)]$, and then compare $X_t$ with Brownian motion with drift.

Let $n \in \mathbb{Z}_+$ be such that $n\delta < 1$, assume $\Delta_n$ holds (note that $\Delta_0$ always holds) and let $t \in [T_{y+n\delta}, U(\delta, T_{y+n\delta}, X)]$. Then

$$g_t(X_t) \geq g_{T_{y+n\delta}}(X_t) - \frac{\delta^2}{2} \geq g_{T_{y+(n-1)+\delta}}(X_t) - \delta^2 \geq g_{T_y}(X_t) - \delta^2, \tag{18}$$

the last inequality following, when $n \geq 1$, from $X_t \geq y + (n-1)^+\delta \geq \sup_{u \in [T_y, T_{y+(n-1)+\delta}]} X_u$.

Now, by inequality (6),

$$g_{T_y}(X_t) \geq \left(\int_{y-\delta}^{y} l_y(v)\,dv\right) \min((f(n-1)\delta), f((n+2)\delta))$$

$$\geq \frac{1}{\delta^2}\left(-\mathbb{1}_{\{n=0\}}\delta + \mathbb{1}_{\{n\neq 0\}}\frac{(n-1)\delta}{10}\right) \geq \frac{1}{\delta}\left(\frac{n}{10} - 1\right), \tag{19}$$

using in the second inequality that $f(x) \geq x/10$ for all $x \in [0,3]$ (since $1 + |x|^{1+\beta} \leq 1 + 3^{1+\beta} \leq 10$). Note that Eq. (19) remains true in the case $\delta = 1$ (which implies $n = 0$ since $n\delta < 1$ by assumption), since $\int_{y-1}^{y} l_y(v)\,dv \leq 1$.

In summary, Eqs (18) and (19) imply, using $\delta \leq 1$,

$$g_t(X_t) \geq \frac{1}{\delta}\left(\frac{n}{10} - 2\right). \tag{20}$$

Given $n \in \mathbb{Z}_+$, let $h := (n/10 - 2)/\delta$ and let $M_t := B_t + ht$ (with $B_0 = 0$). Then a comparison result (cf. [21] for instance) yields that, if $T_{y+n\delta} < \infty$, then

$$\mathbb{P}[T_{y+(n+1)\delta} = U(\delta, T_{y+n\delta}, X)|\mathcal{F}_{T_{y+n\delta}}] \geq \mathbb{P}[M_{U(\delta,0,M)} = \delta]. \tag{21}$$

If $n \geq 80$, Lemma 7 provides (with $a := \delta$ and $h$ defined above)

$$\mathbb{P}[M_{U(\delta,0,M)} = \delta] \geq 1 - \exp\left(4 - \frac{n}{10}\right). \tag{22}$$

If $n < 80$, Girsanov's lemma implies

$$\mathbb{P}[M_{U(\delta,0,M)} = \delta] \geq \exp\left[h\delta - \frac{h^2}{2}\frac{\delta^2}{2}\right]\mathbb{P}[M_{U(\delta,0,B)} = \delta] \geq e^{-3}\frac{1}{2}\mathbb{P}[|B_{\delta^2/2}| \geq \delta] = \text{Cst} \tag{23}$$

the last inequality following from the self-similarity of the Brownian motion.

In summary, Eqs (17), (21)–(23) yield, for all $y \in \mathbb{Z}_+$ such that $T_y < \infty$,

$$\mathbb{P}(\mathcal{E}_y|\mathcal{F}_{T_y}) \geq \text{Cst} \prod_{k \geq 80}\left(1 - \exp\left(4 - \frac{n}{10}\right)\right) \geq \text{Cst} > 0.$$

*4.2. Proof of Lemma 5*

Let, for all $x \in \mathbb{Z}_+$,

$$C_x := A_x \cap \left\{\inf_{t \in [T_{x-1}, T_x]} X_t \geq x - \frac{3}{2}\right\} \cap \left\{\inf_{t \in [T_{x-1}, T_x]} \frac{G(t)}{t^{\alpha-1}} \geq \frac{(4c_0)^{-\beta}}{2}\right\}.$$

We divide the proof in two parts. In part (1) we prove that, for any sufficiently large $x \in \mathbb{Z}_+$,

$$A_x \subseteq C_{x+1}. \tag{24}$$

This implies, since $A_x$ holds infinitely often by Lemma 3, that $C_x$ holds for any large $x$. In part (2) we assume $C_x$ for any sufficiently large $x$, and conclude that $\Gamma_+ \cap \{\liminf G(t)/t^{\alpha-1} > 0\}$ holds.

Part (1). Let, for all $x \in \mathbb{Z}_+$,

$$\mathcal{E}_x := \left\{\inf_{0 \leq s - T_x \leq 3(4c_0)^\beta T_x^{1-\alpha}} (B_s - B_{T_x}) \geq -\frac{1}{2}\right\},$$

with the convention that $\mathcal{E}_x = \Omega$ if $T_x = \infty$. Let us first prove in part (1)(a) that $\mathcal{E}_x$ holds for any sufficiently large $x$ and in part (1)(b) that $A_x \cap \mathcal{E}_x \subseteq C_{x+1}$ for any sufficiently large $x$.

(1)(a) Using the standard exponential inequality (4), for any sufficiently large $x$,

$$\mathbb{P}[\mathcal{E}_x^c | \mathcal{F}_{T_x}] \leq \exp\left(\frac{-(4c_0)^{-\beta} T_x^{\alpha-1}}{24}\right) \leq \exp\left(\frac{-(4c_0)^{1-\alpha^{-1}-\beta} x^{1-\alpha^{-1}}}{24}\right),$$

where we use that $x = X_{T_x} \leq 2c_0 T_x^\alpha$ by Theorem A. Therefore $\sum \mathbb{P}[\mathcal{E}_x^c | \mathcal{F}_{T_x}] < \infty$, which proves the claim by Borel–Cantelli lemma.

(1)(b) Assume $A_x \cap \mathcal{E}_x$ holds, and that $x$ is large. Then, for all $t \in [T_x, T_x + 1]$, as long as $X_t \in [x - 1/2, x + 1]$,

$$G(t) = g_t(X_t) \geq g_{T_x}(X_t) - (t - T_x)\|f\|_\infty \geq \inf_{y \in [x-1/2, x+1]} g_{T_x}(y) - 1,$$

and, by Lemma 4, for any sufficiently large $x$,

$$\frac{G(t)}{t^{\alpha-1}} \geq \frac{\inf_{y \in [x-1/2, x+1]} g_{T_x}(y) - 1}{(T_x + 1)^{\alpha-1}} \geq \frac{(4c_0)^{-\beta}}{2}. \tag{25}$$

Therefore, using that $\mathcal{E}_x$ holds,

$$X_{T_x + 6(4c_0)^\beta T_x^{1-\alpha}} \geq \left[\frac{x-1}{2} + 3(4c_0)^\beta T_x^{1-\alpha} \frac{(4c_0)^{-\beta}}{2} T_x^{\alpha-1}\right] \wedge (x+1) \geq x + 1,$$

which implies that $T_{x+1} \leq T_x + 3(4c_0)^\beta T_x^{1-\alpha} \leq T_x + 1$ for any sufficiently large $x \in \mathbb{Z}_+$ and, using again Eq. (25), that $C_{x+1}$ holds.

Part (2). Assume that $C_x$ holds for any sufficiently large $x$. The asymptotic estimate of $G(t)$ follows immediately. Let us prove the estimate of $\inf_{0 \leq s \leq t}(g_t(X_t) - g_s(X_t))$ as $t$ goes off to infinity: For all large $t \in \mathbb{R}_+$, there exists $x \in \mathbb{Z}_+$ such that $X_t \in [x, x+1)$, and $x$ is large for $t$ large enough, since $\lim_{t \to \infty} X_t = \infty$ by assumption, so that $C_y$ holds for $y \geq x$. Hence, $t \leq T_{X_t} + 2$ since on one hand, $t \leq T_{x+2} \leq T_x + 2$ ($t > T_{x+2}$ would imply $X_t \geq x + 2 - 1/2 = x + 3/2$, which is contradictory), and since on the other hand $T_{X_t} \geq T_x$ by definition of $T_x$. This implies

$$g_t(X_t) - g_s(X_t) = \int_s^t f(X_t - X_u)\,du = \int_s^{T_{X_t} \vee s} f(X_t - X_u)\,du + \int_{T_{X_t} \vee s}^t f(X_t - X_u)\,du$$

$$\geq \int_{T_{X_t} \vee s}^t f(X_t - X_u)\,du \geq -\|f\|_\infty (t - T_{X_t}) \geq -2\|f\|_\infty \geq -2.$$

*4.3. Proof of Lemma* 6

In this section, we explain why the conclusions of Theorem 4 in [18] almost surely hold on $\Gamma_+$, notwithstanding that assumptions $f$ nonnegative and $f(0) > 0$ are not satisfied. The proof in [18] makes use of these assumptions on two occasions. First, Lemma 5.1 and the beginning of the proof of Lemma 5.2 apply the following inequality, referred to as (3.7): there exists $A > 0$ such that, for all $t \geq s \geq 0$,

$$X_t - X_s \geq (t-s)\frac{A^{1/2}}{2} + \inf_{s \leq r \leq t}(B_t - B_r) + \inf_{s \leq r \leq t}(B_r - B_s) - 1.$$

Our condition $\liminf G(t) > 0$ on $\Gamma_+$ implies that there exist a.s. $A > 0$ and $C \in \mathbb{R}_+$ such that $G(u) \geq A$ for all $u \geq C$, so that for all $t \geq s$,

$$\begin{aligned}X_t - X_s &= \int_s^t G(u)\,\mathrm{d}u + B_t - B_s \geq (t-s)A - C(A+C) + B_t - B_s \\ &\geq (t-s)A + \inf_{s \leq r \leq t}(B_t - B_r) + \inf_{s \leq r \leq t}(B_r - B_s) - C(A+C)\end{aligned}$$

so that the inequality continues to hold with this difference that 1 is replaced by a constant $D := C(A+C) > 0$, which does not modify the consequences when $T$ is large enough, depending on $D$.

Second, the assumption $f$ nonnegative also appears in Eq. (5.8) in [18] (in the proof of Lemma 5.2), which makes use of the following inequality, for $v = t - u_T$ (where $u_T$ is a constant depending on $T$):

$$\int_0^t f(X_{Tt} - X_{Ts})\,\mathrm{d}s \geq \int_0^v f(X_{Tt} - X_{Ts})\,\mathrm{d}s.$$

Observe that

$$\int_0^t f(X_{Tt} - X_{Ts})\,\mathrm{d}s = \frac{1}{T}\int_0^{Tt} f(X_{Tt} - X_s)\,\mathrm{d}s = \frac{g_{Tt}(X_{Tt})}{T}$$

and

$$\int_0^v f(X_{Tt} - X_{Ts})\,\mathrm{d}s = \frac{g_{Tv}(X_{Tt})}{T}.$$

Therefore, on $\Gamma_+$, $\liminf_{t \to \infty}(\inf_{0 \leq s \leq t}(g_t(X_t) - g_s(X_t))) > -\infty$ implies

$$\int_0^t f(X_{Tt} - X_{Ts})\,\mathrm{d}s \geq \int_0^v f(X_{Tt} - X_{Ts})\,\mathrm{d}s - \mathrm{O}(T^{-1}).$$

The error term of order $T^{-1}$ turns out to be negligible with respect to the lower bound of $\int_0^v f(X_{Tt} - X_{Ts})\,\mathrm{d}s$ by $\mathrm{Cst}(\alpha,\gamma,c)t^{\alpha-1}/T^{\alpha\beta}$ later in the proof Lemma 5.2 [18] (Eq. (5.12)). Indeed,

$$T^{-1} = T^{-\alpha\beta}T^{\alpha\beta-1} = T^{-\alpha\beta}T^{1-\alpha} = T^{-\alpha\beta}\mathrm{o}(t^{\alpha-1}),$$

if we suppose that $t \geq t_T = c^{-1/\alpha}T^{-\lambda/\alpha} \gg T^{-1}$, as it is indeed assumed in this part of the proof in [18].

**Acknowledgments**

This research project was supported in part by the Swiss National Science Foundation Grants 200021-107425 (T. Mountford) and 200021-1036251/1 (P. Tarrès).
We are grateful to the associate editor and anonymous referee for very useful comments and suggestions.